\documentclass[11pt,a4paper,reqno]{amsart}%
\usepackage{amsthm,amsmath,amsfonts,amssymb,xcolor,amsxtra,dsfont,bm,mathrsfs}
\usepackage[latin1]{inputenc}
\usepackage{color}
\usepackage{amssymb}\usepackage{graphicx}
\usepackage{tikz}
\usepackage{algorithm}
\usepackage{algpseudocode}
\usepackage[title]{appendix}
\usepackage{nicematrix}

\usepackage[colorlinks=true,urlcolor=green!50!black,
citecolor=red!60!black,linkcolor=blue!75!black,linktoc=all,pdfpagelabels, bookmarksnumbered,bookmarksopen]{hyperref}

\usepackage{bookmark}

\usepackage[capitalize,nameinlink,noabbrev]{cleveref}

\theoremstyle{plain}

\theoremstyle{remark}

\allowdisplaybreaks

\title[Spectral Algorithms for 3-Wave Kinetic  and $C_{12}$ Quantum Boltzmann]{Spectral Algorithms for 3-Wave Kinetic  and $C_{12}$ Quantum Boltzmann Equations with General Resonance Manifolds in \(\mathbb{R}^d\)
}

\author[T. T. Le]{Thanh Trung Le}
\address{School of Mathematics and Statistics, University of Economics Ho Chi Minh City, Ho Chi Minh City 700000, Vietnam}
\email{thanhtrungle@ueh.edu.vn}
\thanks{T. T. Le is funded by University of Economics Ho Chi Minh City, Vietnam (UEH)}

\author[M.-B. Tran]{Minh-Binh Tran}
\address{Department of Mathematics, Texas A\&M University, College Station, TX 77843, USA}
\email{minhbinh@tamu.edu}
\thanks{M.-B. T is  funded in part by  a   Humboldt Fellowship,   NSF CAREER  DMS-2303146, and NSF Grants DMS-2204795, DMS-2305523,  DMS-2306379.}

\begin{document}
	\date{\today}

	\begin{abstract}
	Following recent developments in numerical schemes for 3-wave kinetic equations \cite{Banks-Shatah-2025-new,Das-Binh-2025-mumerical,Walton-Binh-2023-numerical,Walton-Binh-2024-deep,Walton-Binh-2025-numerical}, we develop spectral algorithms for multidimensional 3-wave kinetic equations and \(C_{12}\) quantum Boltzmann equations with general polynomial dispersion relations.
 The principal numerical difficulty arises from the resonance constraint, which is supported on a nonlinear manifold in wave-vector space. We approximate the corresponding Dirac distribution by means of a truncated Fourier representation and derive two spectral discretizations of the resulting collision operator. The first method is a direct spectral scheme based on a truncated Fourier approximation of the distribution function and has computational complexity
		\(\mathcal{O}\bigl(L(2N)^{3d}\bigr)\), where \(2N\) is the number of Fourier modes in each coordinate direction and \(L\) is the number of time steps. The second method applies multidimensional fast Fourier transforms directly to the kernel-weighted nonlinear terms, reducing the computational complexity to
		\(\mathcal{O}\bigl(L(2N)^{2d}\log(2N)\bigr)\).
		In two-dimensional tests, the two algorithms agree to nearly machine precision, while the fast method achieves speedups ranging from approximately \(93\) to more than \(2200\) over the resolutions considered. To suppress unresolved high-frequency modes, we introduce a stabilization strategy combining a pre-FFT implementation of the classical \(2/3\)-rule with exponential spectral filtering. Numerical simulations in 2 and 3 dimensions capture the gain--loss dynamics of the \(C_{12}\) quantum Boltzmann equation for both rapidly and algebraically decaying initial data. For the 3-wave kinetic equation with a rapidly growing collision kernel, the computations develop strong oscillatory structures and rapid spectral broadening, providing numerical evidence of an apparent transfer of energy toward high frequencies (energy cascade). The results also demonstrate that the dispersion relation and the spatial dimension significantly influence the transient resonant dynamics and the onset of high-frequency growth.
	\end{abstract}
	
	\maketitle
	\allowdisplaybreaks
	\setcounter{tocdepth}{3}
	\tableofcontents

	\section{Introduction}\label{intro}
	Wave turbulence theory provides a central statistical description of weakly nonlinear dispersive waves evolving over long time scales. The theory has proved relevant across a broad spectrum of physical settings, including inertial-wave interactions in rotating fluids, Alfv\'{e}nic turbulence in the solar wind, wave dynamics in magnetized plasmas, and numerous phenomena arising in plasma and fusion physics. Its conceptual foundations can be traced to the pioneering work of Peierls \cite{Peierls:1993:BRK} and were subsequently developed through the influential contributions of Benney and Saffman \cite{benney1966nonlinear}, Zakharov and Falkovich \cite{zakharov1967weak}, Benney and Newell \cite{benney1969random}, and Hasselmann \cite{hasselmann1962non,hasselmann1974spectral}. This body of work culminated in the formulation of wave kinetic equations, particularly the 3-wave and 4-wave kinetic equations, which govern the resonant transfer of energy among weakly coupled wave modes. In recent years, the mathematical foundations of the theory have advanced substantially through a series of breakthrough results by Deng and Hani, who established rigorous derivations and long-time validity results for wave kinetic equations in several important regimes \cite{deng2019derivation,deng2021propagation, deng2021full,deng2023long,deng2023}. Comprehensive accounts of the physical theory, its principal mechanisms, and its diverse applications may be found in \cite{Nazarenko:2011:WT,Pomeau-Binh-2019-statistical,zakharov2012kolmogorov}.
	
	Following recent progress in the construction of numerical schemes for the radial 
	version of the 3-wave kinetic equations
	\cite{Das-Binh-2025-mumerical, Walton-Binh-2023-numerical, Walton-Binh-2024-deep, Walton-Binh-2025-numerical},
	we develop in this work numerical schemes for the more general multidimensional formulation
	of the equations studied in those references. More precisely, we consider the   kinetic equation
	\cite{Pomeau-Binh-2019-statistical}

	\begin{eqnarray}\label{QB}
		\frac{\partial f}{\partial t} (t,k)
		\	&= &\ \mathbb Q[f] (t,k),\ \ \ f(0, k) \ = \ f_{0}(k),
	\end{eqnarray}
	where the forms of $\mathbb Q$ is  given explicitly below
	\begin{eqnarray}
		\label{C12Discrete}\nonumber
		&& \mathbb Q[f](k)  =      \int_{\mathbb{R}^d}\int_{\mathbb{R}^d}\int_{\mathbb{R}^d}{\mathrm{d}k_1\mathrm{d}k_2\mathrm{d}k_3}\delta(k_1-k_2-k_3) \\\nonumber
		&&\times (\delta(k-k_1)-\delta(k-k_2)-\delta(k-k_3))\delta(\omega_1-\omega_2-\omega_3)\\
		&& \times K^{12}(k_1,k_2,k_3)  \Big[f_2f_3-f_1f_2-f_1f_3-{\mathbf C}f_1\Big],
	\end{eqnarray}
	in which $\omega_i$, $f_i$ stand for $\omega(k_i),f(k_i)$, $k\in\mathbb{R}^d$ is the $3$-dimensional non-zero momentum variable.
	The constant \( {\mathbf C} \) takes the value \( 0 \) or \( 1 \). We distinguish the following cases:
	\begin{itemize}
		\item[(i)] When \( {\mathbf C} = 1 \), \eqref{QB} is the quantum Boltzmann equation describing the density function of excitations outside a condensate \cite{Pomeau-Binh-2019-statistical, Binh-Pomeau-boltzmann,reichl2019kinetic}. In this setting, the collision operator is conventionally identified with the \(C_{12}\) term associated with \(1\leftrightarrow 2\) collision processes.
		\item[(ii)] When \( {\mathbf C} = 0 \), \eqref{QB} is a 3-wave kinetic equation \cite{Soffer-Binh-2019-energy}.
	\end{itemize}
	
	The dispersion relation \( \omega(k) \) is a radial function of $|k|$. In this work, \( \omega(k) \)  is chosen to be of polynomial form
	
	\begin{equation}\label{Poli1}
		\omega(k) = |k|^\alpha, \alpha\ge 1,
	\end{equation}
	\begin{equation}\label{Poli2}
		\omega(k) = c_1|k|^\alpha + c_2|k|^\beta, \alpha,\beta\ge 1, c_1,c_2>0.
		\end{equation}
%
%

The 3-wave and quantum Boltzmann equations \eqref{QB} provide effective descriptions of resonant energy transfer in a wide range of weakly nonlinear wave and quantum systems. Their mathematical structure has therefore been investigated in several distinct physical regimes, including phonon scattering in anharmonic crystal lattices \cite{CraciunBinh,EscobedoBinh,GambaSmithBinh,tran2020reaction}, capillary-wave turbulence \cite{nguyen2017quantum}, internal-wave interactions in stratified oceanic flows \cite{GambaSmithBinh,kim2025wave}, acostic waves \cite{staffilani2026entropy}, beam waves \cite{rumpf2025wave}, and kinetic models for Bose--Einstein condensates \cite{cortes2020system,EPV,escobedo2023linearized1,escobedo2023linearized,ToanBinh,nguyen2017quantum,staffilani2025evolution,staffilani2025finite,staffilani2025formation}. On the numerical side, Banks and Shatah developed a direct discretization method based on a piecewise-polynomial approximation of the resonant manifold and numerical quadrature of the collision integral \cite{Banks-Shatah-2025-new}.

In this work, we develop two spectral algorithms for multidimensional 3-wave kinetic and \(C_{12}\) quantum Boltzmann equations with general polynomial dispersion relations of the type  \ref{Poli1}-\ref{Poli2}. The first follows a direct spectral discretization of the solution, whereas the second applies multidimensional Fourier transforms directly to the kernel-weighted nonlinear terms. Because the collision operator depends on three Fourier indices, its
direct evaluation has computational complexity
\[
\mathcal{O}\bigl(L(2N)^{3d}\bigr),
\]
where $L$ is the number of time steps and $N$ is the number of Fourier modes.

This cost becomes
\(\mathcal{O}\bigl(L(2N)^6\bigr)\) in 2 dimensions and
\(\mathcal{O}\bigl(L(2N)^9\bigr)\) in 3 dimensions, making the
direct method prohibitively expensive at moderate resolutions.

To overcome this limitation, we introduce a fundamentally different
formulation. Rather than expanding each occurrence of the distribution
function inside the collision operator and evaluating the resulting
multi-index sums, we assemble the kernel-weighted nonlinear terms on
the \(2d\)-dimensional pair-variable domain and apply multidimensional
FFTs directly to these terms. The collision operator can then be
recovered from selected Fourier coefficients of the transformed
tensors. This reduces the computational complexity to
\[
\mathcal{O}\bigl(L(2N)^{2d}\log(2N)\bigr).
\]
Consequently, the asymptotic acceleration relative to the direct
spectral method is
\[
\mathcal{O}\left(\frac{(2N)^d}{\log(2N)}\right).
\]
Equivalently, the fast method is nearly linear, up to the logarithmic
FFT factor, in the number \((2N)^{2d}\) of grid values required to
represent the collision terms on the pair-variable domain. This reduction provides an FFT-based acceleration specifically adapted
to collision operators with general resonance manifolds. Its computational significance may be compared with the improvement achieved by fast spectral methods for the classical Boltzmann collision operator
\cite{MouhotPareschi2006,FilbetMouhotPareschi:2006:SBE}
over the direct Fourier--Galerkin spectral method of Pareschi and Russo
\cite{PareschiRusso2000}. The formulation and acceleration mechanism developed here are, however, fundamentally different from those used in the classical Boltzmann setting.

Numerical comparisons in 2 dimensions show that the 2 methods agree to approximately machine precision, while the fast algorithm achieves speedups ranging from approximately \(93\) to more than \(2200\) for the resolutions considered. To control the rapid generation of unresolved modes, we further introduce a stabilization procedure combining a pre-FFT implementation of the classical \(2/3\)-rule with exponential spectral filtering. The resulting method remains computationally feasible for both 2- and 3-dimensional simulations. For the \(C_{12}\) quantum Boltzmann equation, the computations capture the competition between gain and loss mechanisms, preserve the expected approximate radial symmetry, and remain spectrally well resolved for both rapidly and algebraically decaying initial data, as theoretically observed in \cite{Alonso-Gamba-Binh-2016-cauchy} (see also \cite{nguyen2017quantum} for a similar situation). For the 3-wave kinetic equation with a rapidly growing collision kernel, the numerical solutions instead develop strong oscillatory structures together with a rapid broadening of the Fourier spectrum. These computations provide numerical evidence of an apparent transfer of energy toward higher frequencies, consistent with the energy-cascade mechanism established theoretically in \cite{Soffer-Binh-2019-energy,staffilani2025finite,staffilani2026entropy}. They also demonstrate that both the dispersion relation and the spatial dimension substantially influence the transient resonant dynamics and the onset of high-frequency growth. Note that the energy-cascade phenomenon has already been observed numerically for the radial version of the 3-wave kinetic equation; see \cite{Das-Binh-2025-mumerical, Walton-Binh-2023-numerical, Walton-Binh-2024-deep, Walton-Binh-2025-numerical}. The present work confirms that this phenomenon persists in the non-radial numerical schemes. 

\section{The algorithms}\label{algorithms}
\subsection{Derivation of the spectral algorithm}
For any test function \(\varphi\in L^\infty(\mathbb{R}^d)\), the weak formulation of the equation is given by
\cite{Toan-Binh-2019-uniform,Soffer-Binh-2018-dynamics,Soffer-Binh-2019-energy}
\begin{equation}\label{Algo:1}
	\begin{aligned}
		\frac{\mathrm d}{\mathrm dt}
		\int_{\mathbb{R}^d} f(t,k)\varphi(k),\mathrm dk
		={}&
		\int_{\mathbb{R}^{3d}}
		\delta(k_1-k_2-k_3)
		\delta(\omega_1-\omega_2-\omega_3)
		K^{12}(k_1,k_2,k_3)
		\\
		&\quad\times
		\Big[
		f_2f_3-f_1f_2-f_1f_3-\mathbf C f_1
		\Big]
		\Big[
		\varphi(k_1)-\varphi(k_2)-\varphi(k_3)
		\Big]
		\,\mathrm dk_1\,\mathrm dk_2\,\mathrm dk_3
		\\
		={}&
		\int_{\mathbb{R}^{2d}}
		\delta\!\Big(
		\omega(k_2+k_3)-\omega(k_2)-\omega(k_3)
		\Big)
		K^{12}(k_2+k_3,k_2,k_3)
		\\
		&\quad\times
		\Big[
		f(k_2)f(k_3)
		-f(k_2+k_3)f(k_2)
		-f(k_2+k_3)f(k_3)
		-\mathbf C f(k_2+k_3)
		\Big]
		\\
		&\quad\times
		\Big[
		\varphi(k_2+k_3)-\varphi(k_2)-\varphi(k_3)
		\Big]
		\,\mathrm dk_2\,\mathrm dk_3 .
	\end{aligned}
\end{equation}

One of the main challenges in numerically simulating solutions of \eqref{QB} is the approximation
of the resonance manifold \cite{Banks-Shatah-2025-new} represented by the delta function
\(\delta(\omega(k_2+k_3)-\omega(k_2)-\omega(k_3))\).
Our approach is to rewrite this delta function as
\[
\delta(\omega(k_2+k_3)-\omega(k_2)-\omega(k_3))
= \frac{1}{2\pi}\int_{-\infty}^{\infty} \mathrm{d}s\,
e^{is(\omega(k_2+k_3)-\omega(k_2)-\omega(k_3))},
\]
which leads to
\begin{equation}\label{Algo:2}
	\begin{aligned}
		\partial_t \int_{\mathbb{R}^d} \mathrm{d}k\, f(t,k)\varphi(k)
		= {} & \int_{\mathbb{R}^d}\!\!\int_{\mathbb{R}^d}
		\mathrm{d}k_2\,\mathrm{d}k_3\,
		\frac{1}{2\pi}\int_{-\infty}^{\infty} \mathrm{d}s\,
		e^{is(\omega(k_2+k_3)-\omega(k_2)-\omega(k_3))}
		K^{12}(k_2+k_3,k_2,k_3) \\
		& \times \Big[
		f(k_2)f(k_3)
		- f(k_2+k_3)f(k_2)
		- f(k_2+k_3)f(k_3)
		- {\mathbf C} f(k_2+k_3)
		\Big] \\
		& \times \Big[
		\varphi(k_2+k_3)-\varphi(k_2)-\varphi(k_3)
		\Big].
	\end{aligned}
\end{equation}

Similar with the case of the classical Boltzmann equation \cite{MouhotPareschi2006,FilbetMouhotPareschi:2006:SBE,MunafHaackGambaMagin:2014:JCP}, we restrict the domain of $k$ to $[-R,R]^d$ and extend the solution $f(t,k)$ to the whole space. We also restrict the domain of $s$ to $[-M,M]$ and obtain
\begin{equation}\label{Algo:3}
	\begin{aligned}
		\partial_t \int_{[-R,R]^d} \mathrm{d}k\, f(t,k)\varphi(k)
		\approx  &    \int_{[-R,R]^d}\!\!\int_{[-R,R]^d}
		\mathrm{d}k_2\,\mathrm{d}k_3\,
		\frac{1}{2\pi}\int_{[-M,M]} \mathrm{d}s\,
		e^{is(\omega(k_2+k_3)-\omega(k_2)-\omega(k_3))}
		\\
		& \times \Big[
		f(k_2)f(k_3)
		- f(k_2+k_3)f(k_2)
		- f(k_2+k_3)f(k_3)
		-{\mathbf C} f(k_2+k_3)
		\Big] \\
		& \times K^{12}(k_2+k_3,k_2,k_3) \Big[
		\varphi(k_2+k_3)-\varphi(k_2)-\varphi(k_3)
		\Big]\\
		\approx  &    \int_{[-R,R]^d}\!\!\int_{[-R,R]^d}
		\mathrm{d}k_2\,\mathrm{d}k_3\,
		\frac{1}{ \pi}\frac{\sin \big(M(\omega(k_2+k_3)-\omega(k_2)-\omega(k_3))\big)}{\omega(k_2+k_3)-\omega(k_2)-\omega(k_3)}
		\\
		& \times \Big[
		f(k_2)f(k_3)
		- f(k_2+k_3)f(k_2)
		- f(k_2+k_3)f(k_3)
		-{\mathbf C} f(k_2+k_3)
		\Big] \\
		& \times K^{12}(k_2+k_3,k_2,k_3)  \Big[
		\varphi(k_2+k_3)-\varphi(k_2)-\varphi(k_3)
		\Big].
	\end{aligned}
\end{equation}
Setting
\begin{equation}
	\label{Algo:4}
	\begin{aligned}
		\mathscr	W(k_2,k_3) := {} &
		\frac{1}{\pi}
		\frac{\sin\!\big(M(\omega(k_2+k_3)-\omega(k_2)-\omega(k_3))\big)}
		{\omega(k_2+k_3)-\omega(k_2)-\omega(k_3)}
		\,K^{12}(k_2+k_3,k_2,k_3) \\
		= {} &
		\frac{1}{\pi}
		\frac{\sin\!\big(M(\omega(k_2+k_3)-\omega(k_2)-\omega(k_3))\big)}
		{\omega(k_2+k_3)-\omega(k_2)-\omega(k_3)}
		\Big|\omega(k_2+k_3)\,\omega(k_2)\,\omega(k_3)\Big|^{\rho},
	\end{aligned}
\end{equation}
we obtain
\begin{equation}\label{Algo:5}
	\begin{aligned}
		\partial_t \int_{[-R,R]^d} \mathrm{d}k\, f(t,k)\varphi(k)
		\approx {} &
		\int_{[-R,R]^d}\!\!\int_{[-R,R]^d}
		\mathrm{d}k_2\,\mathrm{d}k_3\,
		\mathscr	W(k_2,k_3) \\
		& \times \Big[
		f(k_2)f(k_3)
		- f(k_2+k_3)f(k_2)
		- f(k_2+k_3)f(k_3)
		- {\mathbf C}f(k_2+k_3)
		\Big] \\
		& \times \Big[
		\varphi(k_2+k_3)-\varphi(k_2)-\varphi(k_3)
		\Big].
	\end{aligned}
\end{equation}

We next approximate \(f\) by a truncated Fourier series,
\begin{equation}
	\label{Algo:6}
	f(t,k) \approx f_N(t,k)
	= \sum_{n=-N}^{N-1} {\bf a}_n(t) \, e^{i\frac{\pi}{R} n\cdot k},
	\qquad
	{\bf a}_n(t)
	= \frac{1}{(2R)^d}\int_{[-R,R]^d} \mathrm{d}k \,
	f(t,k)\, e^{-i\frac{\pi}{R} n\cdot k}, \ \ t > 0,
\end{equation}
where the sum is taken over the lattice
\begin{equation}	\label{Algo:7}
	\big\{ n=(n^1,\cdots,n^d)\in\mathbb{Z}^d \;:\; -N \le n^1,\cdots,n^d \le N-1 \big\}.
\end{equation}

Plugging \eqref{Algo:6} into the left-hand side of \eqref{Algo:5} and choosing \(\varphi(k)=e^{-i\frac{\pi}{R} n\cdot k}\), we obtain the following system:
\begin{equation}\label{Algo:ODEs}
	\begin{aligned}
		(2R)^d \dot{{\bf a}}_n(t) \approx {} &
		\int_{[-R,R]^d}\!\!\int_{[-R,R]^d}
		\mathrm{d}k_2\,\mathrm{d}k_3\,
		\mathscr	W(k_2,k_3) \\
		& \times \Big[
		f(k_2)f(k_3)
		- f(k_2+k_3)f(k_2)
		- f(k_2+k_3)f(k_3)
		- {\mathbf C}f(k_2+k_3)
		\Big] \\
		& \times \Big[
		e^{-i\frac{\pi}{R} n\cdot (k_2 + k_3)} - e^{-i\frac{\pi}{R} n \cdot k_2}-e^{-i\frac{\pi}{R} n\cdot k_3} \Big]. \\
	\end{aligned}
\end{equation}

To solve the ODE system \eqref{Algo:ODEs}, we develop two distinct approaches:
\begin{itemize}
	\item The first is a direct spectral method in which the function \(f\) on the right-hand side of \eqref{Algo:ODEs} is approximated by a truncated Fourier series; see \cref{sec:algorithm_1}. This type of approximation is standard in spectral methods for the classical Boltzmann equation; see, for example, \cite{PareschiRusso2000,GambaTharkabhushanam2009,MouhotPareschi2006,FilbetMouhotPareschi:2006:SBE,MunafHaackGambaMagin:2014:JCP}.

	\item The second approach follows a fundamentally different principle. Rather than approximating \(f\) itself by a truncated Fourier series, we apply the Fourier transform directly to the kernel-weighted nonlinear terms, namely, to the products of \(\mathscr W\) with the various combinations of \(f\) appearing on the right-hand side of \eqref{Algo:ODEs}. To the best of our knowledge, this formulation has not previously been investigated in the literature. The resulting algorithm is developed in \cref{sec:algorithm_2}.

\end{itemize}

\subsection{Algorithm 1: a direct spectral method} \label{sec:algorithm_1}
Plugging the truncated Fourier series approximation of $f$ (see \eqref{Algo:6}) into the right-hand side of \eqref{Algo:ODEs}, we obtain the following system
\begin{equation}\label{Algo:8}
	\begin{aligned}
		\dot{{\bf a}}_n(t)
		\approx {} &
		\frac{1}{(2R)^d}\int_{[-R,R]^d}\!\!\int_{[-R,R]^d}
		\mathrm{d}k_2\,\mathrm{d}k_3\,
		\mathscr{W}(k_2,k_3) \\
		& \times \Bigg[
		\sum_{p=-N}^{N-1} {\bf a}_p(t)\, e^{i\frac{\pi}{R} p\cdot k_2}
		\sum_{m=-N}^{N-1} {\bf a}_m(t) \, e^{i\frac{\pi}{R} m\cdot k_3}
		- \sum_{p=-N}^{N-1} {\bf a}_p(t) \, e^{i\frac{\pi}{R} p\cdot (k_2+k_3)}
		\sum_{m=-N}^{N-1} {\bf a}_m(t) \, e^{i\frac{\pi}{R} m\cdot k_2} \\
		& \qquad
		- \sum_{p=-N}^{N-1} {\bf a}_p(t) \, e^{i\frac{\pi}{R} p\cdot (k_2+k_3)}
		\sum_{m=-N}^{N-1} {\bf a}_m(t) \, e^{i\frac{\pi}{R} m\cdot k_3}
		-{\mathbf C} \sum_{p=-N}^{N-1} {\bf a}_p(t) \, e^{i\frac{\pi}{R} p\cdot (k_2+k_3)}
		\Bigg] \\
		& \times \Big[
		e^{-i\frac{\pi}{R} n\cdot (k_2+k_3)}
		- e^{-i\frac{\pi}{R} n\cdot k_2}
		- e^{-i\frac{\pi}{R} n\cdot k_3}
		\Big].
	\end{aligned}
\end{equation}
Replacing $\approx$ by $=$, and developing \eqref{Algo:8}, we find a system of ODEs with unknowns $\{{\bf a}_n\}_{n=-N}^{N-1}$
\begin{equation}\label{Algo:9}
	\begin{aligned}
		\dot{{\bf a}}_n
		= {} &
		\frac{1}{(2R)^d}\sum_{p=-N}^{N-1} 	\sum_{m=-N}^{N-1} {\bf a}_p{\bf a}_m\Bigg[	\int_{[-R,R]^d}\!\!\int_{[-R,R]^d}
		\mathrm{d}k_2\,\mathrm{d}k_3\,
		\mathscr{W}(k_2,k_3) \Bigg] \\
		& \times \Bigg[
		e^{i\frac{\pi}{R} p\cdot k_2}
		e^{i\frac{\pi}{R} m\cdot k_3}
		-   e^{i\frac{\pi}{R} p\cdot (k_2+k_3)}
		e^{i\frac{\pi}{R} m\cdot k_2}  -    e^{i\frac{\pi}{R} p\cdot (k_2+k_3)}  e^{i\frac{\pi}{R} m\cdot k_3}
		\Bigg] \\
		& \times \Big[
		e^{-i\frac{\pi}{R} n\cdot (k_2+k_3)}
		- e^{-i\frac{\pi}{R} n\cdot k_2}
		- e^{-i\frac{\pi}{R} n\cdot k_3}
		\Big]\\
		&- \frac{1}{(2R)^d}	{\mathbf C}\sum_{p=-N}^{N-1} 	  {\bf a}_p \Bigg[	\int_{[-R,R]^d}\!\!\int_{[-R,R]^d}
		\mathrm{d}k_2\,\mathrm{d}k_3\,
		\mathscr{W}(k_2,k_3)  e^{i\frac{\pi}{R} p\cdot (k_2+k_3)}
		\Bigg]\\	& \times \Big[
		e^{-i\frac{\pi}{R} n\cdot (k_2+k_3)}
		- e^{-i\frac{\pi}{R} n\cdot k_2}
		- e^{-i\frac{\pi}{R} n\cdot k_3}
		\Big]\Bigg]  .
	\end{aligned}
\end{equation}
We reduce \eqref{Algo:9} to
\begin{equation}\label{Algo:11}
	\begin{aligned}
		\dot{\mathbf{a}}_n
		= {} &
		\sum_{p=-N}^{N-1}\sum_{m=-N}^{N-1}
		\mathbf{a}_p\,\mathbf{a}_m\,\mathscr{K}_{n,m,p}
		+ {\mathbf C}\sum_{p=-N}^{N-1}
		\mathbf{a}_p\,\mathscr{H}_{n,p},
	\end{aligned}
\end{equation}
by setting
\begin{equation}\label{Algo:K_nmp}
	\begin{aligned}
		\mathscr{K}_{n,m,p}
		:= {} &
		\frac{1}{(2R)^d}\int_{[-R,R]^d}\!\!\int_{[-R,R]^d}
		\mathrm{d}k_2\,\mathrm{d}k_3\,
		\mathscr{W}(k_2,k_3) \\
		& \times \Bigg[
		e^{i\frac{\pi}{R} p\cdot k_2}
		e^{i\frac{\pi}{R} m\cdot k_3}
		- e^{i\frac{\pi}{R} p\cdot (k_2+k_3)}
		e^{i\frac{\pi}{R} m\cdot k_2}
		- e^{i\frac{\pi}{R} p\cdot (k_2+k_3)}
		e^{i\frac{\pi}{R} m\cdot k_3}
		\Bigg] \\
		& \times \Bigg[
		e^{-i\frac{\pi}{R} n\cdot (k_2+k_3)}
		- e^{-i\frac{\pi}{R} n\cdot k_2}
		- e^{-i\frac{\pi}{R} n\cdot k_3}
		\Bigg] \\
		=  &
		\frac{1}{(2R)^d}\int_{[-R,R]^d}\!\!\int_{[-R,R]^d}
		\mathrm{d}k_2\,\mathrm{d}k_3\, \mathscr{W}(k_2,k_3) \\
		&\times\Big[
		\;e^{-i\frac{\pi}{R}(n-p)\cdot k_2}\,
		e^{-i\frac{\pi}{R}(n-m)\cdot k_3}
		- e^{-i\frac{\pi}{R}(n-p)\cdot k_2}\,
		e^{-i\frac{\pi}{R}(-m)\cdot k_3}
		- e^{-i\frac{\pi}{R}(-p)\cdot k_2}\,
		e^{-i\frac{\pi}{R}(n-m)\cdot k_3}
		\\
		&\quad
		- e^{-i\frac{\pi}{R}(n-p-m)\cdot k_2}\,
		e^{-i\frac{\pi}{R}(n-p)\cdot k_3}
		+ e^{-i\frac{\pi}{R}(n-p-m)\cdot k_2}\,
		e^{-i\frac{\pi}{R}(-p)\cdot k_3}
		+ e^{-i\frac{\pi}{R}(-p-m)\cdot k_2}\,
		e^{-i\frac{\pi}{R}(n-p)\cdot k_3}
		\\
		&\quad
		- e^{-i\frac{\pi}{R}(n-p)\cdot k_2}\,
		e^{-i\frac{\pi}{R}(n-p-m)\cdot k_3}
		+ e^{-i\frac{\pi}{R}(n-p)\cdot k_2}\,
		e^{-i\frac{\pi}{R}(-p-m)\cdot k_3}
		+ e^{-i\frac{\pi}{R}(-p)\cdot k_2}\,
		e^{-i\frac{\pi}{R}(n-p-m)\cdot k_3}
		\;\Big] \\
	\end{aligned}
\end{equation}
and
\begin{equation}\label{Algo:H_np}
	\begin{aligned}
		\mathscr{H}_{n,p}
		:= &
		-\frac{1}{(2R)^d} \int_{[-R,R]^d}\!\!\int_{[-R,R]^d}
		\mathrm{d}k_2\,\mathrm{d}k_3\,
		\mathscr{W}(k_2,k_3)\,
		e^{i\frac{\pi}{R} p\cdot (k_2+k_3)} \\
		& \times \Big[
		e^{-i\frac{\pi}{R} n\cdot (k_2+k_3)}
		- e^{-i\frac{\pi}{R} n\cdot k_2}
		- e^{-i\frac{\pi}{R} n\cdot k_3}
		\Big] \\
		= & -\frac{1}{(2R)^d} \int_{[-R,R]^d}\!\!\int_{[-R,R]^d}
		\mathrm{d}k_2\,\mathrm{d}k_3\,
		\mathscr{W}(k_2,k_3)\, \\
		&\times \Big[ e^{-i\frac{\pi}{R}(n-p)\cdot k_2}\,
		e^{-i\frac{\pi}{R}(n-p)\cdot k_3}
		- e^{-i\frac{\pi}{R}(n-p)\cdot k_2}\,
		e^{-i\frac{\pi}{R}(-p)\cdot k_3}
		- e^{-i\frac{\pi}{R}(-p)\cdot k_2}\,
		e^{-i\frac{\pi}{R}(n-p)\cdot k_3} \Big].
	\end{aligned}
\end{equation}

Now, we approximate \(W\) by a truncated Fourier series,
\begin{equation}
	W(k_2,k_3)=\sum_{\alpha= -N}^{N-1} \sum_{\beta= -N}^{N-1}
	\widehat W_{\alpha,\beta}\,
	e^{i\frac{\pi}{R}\alpha\cdot k_2}\,
	e^{i\frac{\pi}{R}\beta\cdot k_3}.
\end{equation}
where
\begin{equation} \label{eqn:Fourier_W}
	\widehat W (\alpha,\beta) =
	\frac{1}{(2R)^{2d}}
	\int_{[-R,R]^d}\!\!\int_{[-R,R]^d}
	W(k_2,k_3)\,
	e^{-i\frac{\pi}{R}(\alpha\cdot k_2+\beta\cdot k_3)}
	\,dk_2\,dk_3.
\end{equation}
Hence, we deduce from \eqref{Algo:K_nmp} that
\begin{equation} \label{eqn:K_FW}
	\begin{aligned}
		\mathscr{K}_{n,m,p}	
		= & (2R)^d \Big[ \widehat W (n - p, n - m) - \widehat W (n - p, -m) - \widehat W (-p, n - m) \\
		&\qquad - \widehat W (n - p -m, n - p) + \widehat W (n - p -m, -p) + \widehat W (- p - m, n - p) \\
		&\qquad - \widehat W (n - p, n - p - m) + \widehat W (n - p, - p - m) + \widehat W (-p, n - p - m) \Big].
	\end{aligned}
\end{equation}
Similarly, we obtain from \eqref{Algo:H_np} that
\begin{equation} \label{eqn:H_FW}
	\begin{aligned}
		\mathscr{H}_{n,p}
		=  - (2R)^d \Big[ \widehat W (n - p, n - p) - \widehat W (n - p, -p) - \widehat W (-p, n-p)  \Big].
	\end{aligned}
\end{equation}
Note that in some cases, the indices of $\widehat{W}$ in \eqref{eqn:K_FW} and \eqref{eqn:H_FW} may fall outside the range $[-N,\,N-1]$. In such cases, we simply take the indices modulo $2N$.

At this point, the ODE system \eqref{Algo:11} can be integrated in time using a forward Euler scheme. We fix a time step \(\Delta t\) and partition the time interval \([0,T]\) into \(L\) subintervals, with \(T = L\Delta t\). 
For \(j = 0,1,2,\ldots L - 1\), we compute
\begin{equation}\label{Algo:12}
	\begin{aligned}
		\mathbf{a}_n((j+1)\Delta t)
		= {} & \mathbf{a}_n(j\Delta t)
		+ \Delta t \sum_{p=-N}^{N-1}\sum_{m=-N}^{N-1}
		\mathbf{a}_p(j\Delta t)\,\mathbf{a}_m(j\Delta t)\,\mathscr{K}_{n,m,p} \\
		& \quad
		+{\mathbf C} \Delta t \sum_{p=-N}^{N-1}
		\mathbf{a}_p(j\Delta t)\,\mathscr{H}_{n,p},
	\end{aligned}
\end{equation}
where $\mathscr{K}_{n,m,p}$ and $\mathscr{H}_{n,p}$ are respectively evaluated by \eqref{eqn:K_FW} and \eqref{eqn:H_FW}, with initial datum $f_0$ is also approximated by the truncated Fourier series as follows
\begin{equation}
	\label{Algo:initial_condition}
	f_0(k) \approx f_N^0(t,k)
	= \sum_{n=-N}^{N-1} {\bf a}_n^0(t) \, e^{i\frac{\pi}{R} n\cdot k},
	\qquad
	{\bf a}_n^0(t)
	= \frac{1}{(2R)^d}\int_{[-R,R]^d} \mathrm{d}k \,
	f_0(k)\, e^{-i\frac{\pi}{R} n\cdot k}, \ \ t > 0.
\end{equation}
The detailed algorithm is presented in \cref{algorithm_1}.

\begin{algorithm}
	\caption{The direct spectral algorithm}
	\label{algorithm_1}
	\begin{algorithmic}[1]
		\State Precompute the value of $\widehat W$ as defined in \eqref{eqn:Fourier_W} using the fast Fourier transform (FFT) - cost $\mathcal{O}((2N)^{2d} \log(2N))$.
		\State Precompute the  initial condition \(\mathbf{a}_n^0\), obtained from the initial datum \(f_0\) via \eqref{Algo:initial_condition} using FFT - cost $\mathcal{O}((2N)^d \log (2N))$.
		\State For each time step \(j = 0,1,2,\ldots L - 1\), using a forward Euler scheme \eqref{Algo:12}, compute $\mathbf{a}_n((j+1)\Delta t)$ - cost $\mathcal{O}(L (2N)^{3d})$.
		\State For each time step, reconstruct
		\(f((j+1)\Delta t,\cdot)\) from \(\{\mathbf{a}_n((j+1)\Delta t)\}\) using FFT  - cost $\mathcal{O}(L(2N)^d \log (2N))$.
	\end{algorithmic}
\end{algorithm}

\begin{itemize}
	\item \textbf{Complexity.} Step~3 is by far the most computationally expensive step. Therefore, the overall computational complexity of \cref{algorithm_1} is $\mathcal{O}(L (2N)^{3d})$. For instance, this complexity becomes $\mathcal{O}(L(2N)^6)$ in two dimensions and $\mathcal{O}(L(2N)^9)$ in three dimensions.
	\item \textbf{Memory-efficient implementation.} In Step~1, rather than precomputing the tensors
	$\mathscr{K}_{n,m,p}$ and $\mathscr{H}_{n,p}$ directly, we first
	compute $\widehat{W}$ using the fast Fourier transform (FFT) as
	defined in \eqref{eqn:Fourier_W}. The tensors
	$\mathscr{K}_{n,m,p}$ and $\mathscr{H}_{n,p}$ are then evaluated
	using \eqref{eqn:K_FW} and \eqref{eqn:H_FW}. This approach significantly reduces both the computational cost and the memory requirement.
	For example, in three dimensions, the precomputation of the tensor $\mathscr{K}_{n,m,p}$ requires $\mathcal{O}((2N)^9)$ operations. More importantly, storing the $9$-dimensional tensor $\mathscr{K}_{n,m,p}$ requires approximately $1.1$ terabytes (TiB) of memory when $2N=16$ and the data type is \texttt{std::complex<double>}. Such a memory requirement exceeds the capacity of many modern supercomputers. In contrast, storing the $6$-dimensional tensor $\widehat{W}$ requires only approximately $256$~MiB when $2N=16$ and $16$~GiB when $2N=32$, making the proposed implementation feasible for large-scale three-dimensional computations.
\end{itemize}

\subsection{Algorithm 2: a fast spectral method} \label{sec:algorithm_2}
Unlike the direct spectral method presented in
\cref{sec:algorithm_1}, we apply the Fourier transform directly to the kernel-weighted terms. To this end, we introduce the following auxiliary functions:
\begin{equation} \label{eqn:intermediate_tensors}
	\begin{aligned}
		T(k_2, k_3) &= \mathscr	W(k_2,k_3) f(k_2)f(k_3), \\
		Y(k_2, k_3) &= \mathscr	W(k_2,k_3) f(k_2 + k_3)f(k_2), \\
		U(k_2, k_3) &= \mathscr	W(k_2,k_3) f(k_2 + k_3)f(k_3), \\
		I(k_2, k_3) &= \mathscr	W(k_2,k_3) f(k_2 + k_3),
	\end{aligned}
\end{equation}
Let $\widehat T$ be a Fourier transform of $T$, i.e.,
\begin{equation}
	\widehat T (\alpha,\beta) =
	\frac{1}{(2R)^{2d}}
	\int_{[-R,R]^d}\!\!\int_{[-R,R]^d} dk_2\,dk_3\,
	T(k_2,k_3)\,
	e^{-i\frac{\pi}{R}(\alpha\cdot k_2+\beta\cdot k_3)}.
\end{equation}
Similarly,  $\widehat Y, \widehat U, \widehat I$ are the Fourier transform of $Y, U, I$. The ODE system \eqref{Algo:ODEs} can be rewrite as follows,
\begin{equation}
	\begin{aligned}
		\dot{{\bf a}}_n(t) \approx {} (2R)^{d} \Bigg \{ &\Big[  \widehat T (n, n) - \widehat T (n, 0) - \widehat T (0, n) \Big] 
		-  \Big[  \widehat Y (n, n) - \widehat Y (n, 0) - \widehat Y (0, n) \Big] \\
		- &  \Big[  \widehat U (n, n) - \widehat U (n, 0) - \widehat U (0, n) \Big] 
		-  {\mathbf C}  \Big[  \widehat I (n, n) - \widehat I (n, 0) - \widehat I (0, n) \Big] \Bigg \}.
	\end{aligned}
\end{equation}
Using forward Euler scheme with initial condition \eqref{Algo:initial_condition} , for \(j = 0,1,2,\ldots L - 1\), we get that
\begin{equation}
	\begin{aligned} \label{Algo:15}
		\mathbf{a}_n((j+1)\Delta t)
		= {}  \mathbf{a}_n(j\Delta t)
		+ (2R)^d \Delta t \Bigg \{ &\Big[  \widehat T (n, n) - \widehat T (n, 0) - \widehat T (0, n) \Big] - \Big[  \widehat Y (n, n) - \widehat Y (n, 0) - \widehat Y (0, n) \Big] \\
		-  &\Big[  \widehat U (n, n) - \widehat U (n, 0) - \widehat U (0, n) \Big] -  {\mathbf C} \Big[  \widehat I (n, n) - \widehat I (n, 0) - \widehat I (0, n) \Big] \Bigg\}
	\end{aligned}
\end{equation}
The detailed algorithm is presented in \cref{algorithm_2}.

\begin{algorithm}
	\caption{The fast spectral algorithm}
	\label{algorithm_2}
	\begin{algorithmic}[1]
		\State Precompute the value of tensor $\mathscr	W$ using \eqref{Algo:4} - cost $\mathcal{O}((2N)^{2d}$.
		\State Precompute the  initial condition \(\mathbf{a}_n^0\), obtained from the initial datum \(f_0\) via \eqref{Algo:initial_condition} using FFT - cost $\mathcal{O}((2N)^d \log (2N))$.
		\State For each time step \(j = 0,1,2,\ldots L - 1\), construct the tensors $T, Y, U, I$ via the value of $\mathscr W$ and $f(j\Delta t)$ - cost $\mathcal{O}(L(2N)^{2d})$.
		\State For each time step, compute $\widehat T$, $\widehat Y$, $\widehat U$, $\widehat I$ using FFT - cost $\mathcal{O}(L (2N)^{2d} \log(2N))$.
		\State For each time step, compute $\mathbf{a}_n((j+1)\Delta t)$ using \eqref{Algo:15} - cost $\mathcal{O}(L(2N)^{d})$.
		\State For each time step, reconstruct
		\(f((j+1)\Delta t,\cdot)\) from \(\{\mathbf{a}_n((j+1)\Delta t)\}\) using FFT  - cost $\mathcal{O}(L(2N)^d \log (2N))$.
	\end{algorithmic}
\end{algorithm}

\begin{itemize}
	\item \textbf{Complexity.} Step~4 is by far the most computationally expensive step. Therefore, the overall computational complexity of \cref{algorithm_2} is $\mathcal{O}(L (2N)^{2d} \log(2N))$. Compared with \cref{algorithm_1}, the \cref{algorithm_2} reduces the overall computational complexity from $\mathcal{O}(L(2N)^{3d})$ to $\mathcal{O}(L(2N)^{2d}\log(2N))$. The computational gain achieved by the proposed method is comparable in practical significance to the improvement obtained, for the classical Boltzmann collision operator, by replacing the direct Fourier--Galerkin method of Pareschi and Russo \cite{PareschiRusso2000} with the fast spectral methods developed in \cite{MouhotPareschi2006,FilbetMouhotPareschi:2006:SBE}. This comparison concerns only the impact of the complexity reduction: the formulation introduced here and the mechanism underlying its FFT-based acceleration are fundamentally different from those used in the classical Boltzmann setting.
	
	\item \textbf{Memory requirement.} In the three-dimensional case, the \cref{algorithm_2} requires two six-dimensional tensors in memory: one for storing $\mathscr W$ in Step~1 and one reusable working tensor for the functions $T$, $Y$, $U$, and $I$, together with their Fourier transforms. Consequently, the memory requirement is approximately twice that of \cref{algorithm_1}. Nevertheless, this increase remains practical; for example, the required memory is about $512$~MiB for $2N=16$ and $32$~GiB for $2N=32$, which is well within the capabilities of modern computing platforms.
\end{itemize}

\section{Numerical Results}
This section presents a comprehensive evaluation of the proposed spectral algorithms through a series of numerical experiments. \Cref{sec:setup} first introduces the common numerical settings used throughout all simulations.
\Cref{sec:performance} compares the computational performance of \cref{algorithm_1} and \cref{algorithm_2} under identical conditions, without any stabilization techniques, in order to assess their intrinsic efficiency.

Based on the performance comparison in \cref{sec:performance},
\cref{algorithm_2} is adopted in the remainder of this paper.
\Cref{sec:stabilization} introduces the spectral stabilization
strategy adopted in \cref{algorithm_2}, including the pre-FFT application of the classical $2/3$-rule and the exponential spectral filter.
Finally, \cref{sec:qBe} and \cref{sec:3wave} demonstrate the
performance and robustness of the stabilized \cref{algorithm_2} through simulations of the quantum Boltzmann equation and the three-wave kinetic equation, respectively.

\subsection{Numerical setup} \label{sec:setup}

As discussed above, and in analogy with spectral treatments of the classical Boltzmann equation \cite{MouhotPareschi2006,FilbetMouhotPareschi:2006:SBE,MunafHaackGambaMagin:2014:JCP}, we truncate the wave-vector domain to the computational box $[-R,R]^d$. Following the numerical settings adopted in these works, we take $R=10$, which is sufficiently large for the simulations considered here. We extend the truncated solution $f(t,k)$ to the whole space as described above. Similarly, we restrict the auxiliary variable $s$ to the interval $[-M,M]$. In the numerical experiments presented below, we set $M=4$. These particular values are chosen for convenience; the implementation applies equally well to arbitrary choices of the truncation parameters $R$ and $M$.

For all algorithms, time integration is performed using the explicit Euler method with a fixed time step $\Delta t=0.005$.

The algorithms are implemented in C++ using double-precision arithmetic. All computations are performed on a workstation equipped with an AMD Ryzen 9 7950X 16-Core Processor. To reduce the wall-clock time, shared-memory parallelization based on OpenMP is employed, while all post-processing and visualization are carried out using Python.

Throughout this section, the upper row of each figure displays the logarithmic Fourier spectrum, whereas the lower row shows the corresponding numerical solution in physical space at selected times. The Fourier spectrum is visualized by plotting $\ln(|a_n|+10^{-14})$, where the small positive constant $10^{-14}$ is introduced solely to avoid the logarithmic singularity at vanishing Fourier coefficients and has no influence on the numerical results.

\subsection{Performance Comparison} \label{sec:performance}
This section compares the performance of the proposed \cref{algorithm_2} with \cref{algorithm_1}, which is based on the conventional computational framework widely adopted for the classical Boltzmann equations  \cite{MouhotPareschi2006,FilbetMouhotPareschi:2006:SBE,MunafHaackGambaMagin:2014:JCP}. The comparison focuses on both numerical accuracy and computational efficiency in the 2-dimensional case.

The numerical results first demonstrate that \cref{algorithm_2} preserves the accuracy of \cref{algorithm_1}. As shown in \cref{tab:L2-difference}, the $L_2$-norm of the difference between the 2 numerical solutions remains on the order of $10^{-12}$ for all tested Fourier transform sizes ($2N = 16;\, 32; \, 64; \, 128$), which is close to machine precision. This confirms that the new algorithm produces virtually identical solutions while employing a fundamentally different computational strategy.

\begin{table}[!ht]
	\caption{$L_2$ - norm of the difference of \cref{algorithm_1} and \cref{algorithm_2} in the two-dimensional case.}
	\label{tab:L2-difference}
	\begin{tabular}{|c|c|c|c|c|}
		\hline
		& $2N = 16$ & $2N = 32$ & $2N = 64$ & $2N = 128$ \\ 
		\hline
		$L_2$ - norm & $2.2335 \, e-12$ & $2.4572 \, e-12$ & $2.7848 \, e-12$ & $2.6285 \, e-12$ \\
		\hline
	\end{tabular}
\end{table}

For a fair performance comparison, the execution time is measured as the average CPU time per time step over $100$ consecutive time steps. The results are summarized in \cref{tab:CPU-time}. The measured speedup is defined as the ratio between the average CPU time per time step of \cref{algorithm_1} and that of \cref{algorithm_2}. As the Fourier transform size ($2N$) increases, the computational advantage of \cref{algorithm_2} becomes increasingly significant. The measured speedup increases from approximately 93 times for $2N=16$ to 240 times for $2N=32$ and more than 2200 times for $2N=64$. For the largest tested Fourier transform size ($2N=128$), the execution time of \cref{algorithm_1} becomes prohibitively expensive and is therefore omitted, whereas \cref{algorithm_2} still requires only $19.53$ seconds per time step. The measured performance follows the trend predicted by the complexity analysis presented in \cref{sec:algorithm_1} and \cref{sec:algorithm_2}. Since \cref{algorithm_1} and \cref{algorithm_2} have computational complexities of $\mathcal{O}((2N)^{3d})$ and $\mathcal{O}((2N) ^ {2d} \log (2N))$, respectively, the asymptotic speedup factor scales as
$\mathcal{O}\left(\frac{(2N)^d}{\log(2N)}\right)$.
Accordingly, the asymptotic speedup factor increases from 64 for $2N=16$ to approximately 2340 for $2N=128$, following the same trend as the measured speedup.

\begin{figure}[!ht]
	\centering
	\includegraphics[width=0.8\linewidth]{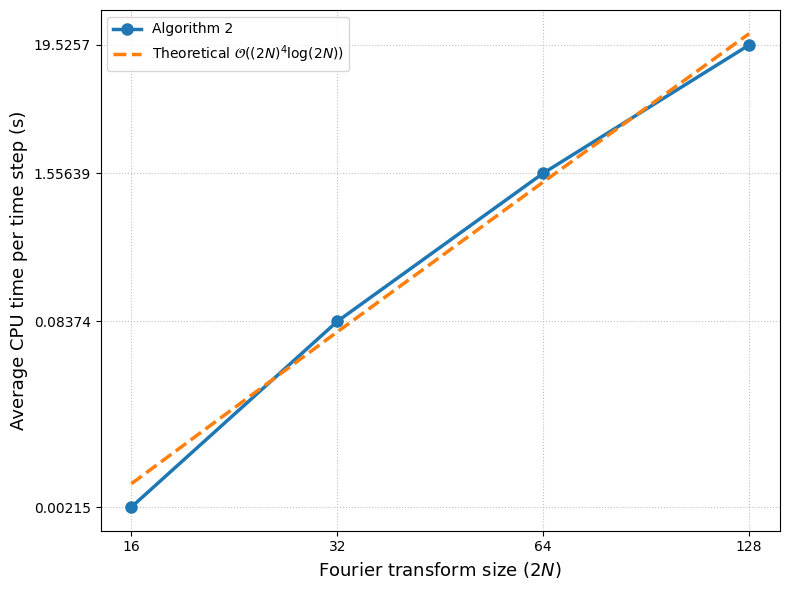}
	\caption{Comparison between the measured CPU time per time step of \cref{algorithm_2} and the theoretical complexity $\mathcal{O}((2N)^4\log(2N))$ for the two-dimensional case. Both axes are shown on logarithmic scales.}
	\label{fig:cpu_time}
\end{figure}

The remaining discrepancy between the measured and asymptotic speedup factors can be attributed to the constant factors hidden in the complexity estimates and to the different computational characteristics of the two algorithms. In \cref{algorithm_1}, the dominant computational cost is incurred in Step~3, which consists of six nested loops evaluating the collision operator directly in two-dimensional case. This step offers little opportunity for further optimization, and its execution time grows rapidly with the Fourier transform size. In contrast, the most expensive part of \cref{algorithm_2} is the computation of multidimensional Fourier transforms of the four-dimensional tensors $T, Y, U$, and $I$. These operations are carried out using highly optimized FFT routines, whose efficiency improves with increasing transform size through optimized execution plans, better cache utilization, and architecture-specific optimizations. Consequently, the measured speedups are consistently larger than those predicted by the asymptotic complexity analysis, particularly for larger Fourier transform sizes. This behavior is further illustrated in \cref{fig:cpu_time}, where the measured execution times closely follow the predicted asymptotic growth. In particular, the measured growth rate is slightly steeper than the asymptotic trend for small Fourier transform sizes and slightly flatter for larger Fourier transform sizes, reflecting the improved practical efficiency of the FFT implementation as the transform size increases.

\begin{table}[!ht]
	\caption{Average CPU time per time step, measured over 100 time steps, for \cref{algorithm_1} and \cref{algorithm_2} in the two-dimensional case.}
	\label{tab:CPU-time}
	\begin{tabular}{|c|c|c|c|c|}
		\hline
		& \cref{algorithm_1} & \cref{algorithm_2} & Measured speedup & Asymptotic speedup factor  \\ \hline
		2N = 16 & 0.200402 s &  0.002149 s & 93.25 & 64 \\ \hline
		2N = 32 & 20.11392 s &  0.083744 s & 240.18 & 204.8 \\ \hline
		2N = 64 & 3558.164 s &  1.556388 s & 2286.17 & 682.67 \\ \hline
		2N = 128 & ---  & 19.52570 s & --- & 2340.57 \\ \hline
	\end{tabular}
\end{table}

Overall, the numerical experiments demonstrate that \cref{algorithm_2} preserves the numerical accuracy of \cref{algorithm_1} while reducing the computational complexity from $\mathcal{O}((2N)^{3d})$ to $\mathcal{O}((2N) ^{2d} \log (2N))$, resulting in substantial computational savings and excellent scalability for large-scale simulations. To further illustrate the scalability of the proposed algorithm, we also performed experiments in the three-dimensional case. The average CPU time per time step of \cref{algorithm_2} is $1.51912$ s for $2N=16$ and $108.094$ s for $2N=32$. Although no direct comparison with \cref{algorithm_1} is provided due to its prohibitive computational cost in three dimensions, these results demonstrate that \cref{algorithm_2} remains computationally feasible for higher-dimensional problems.


\subsection{Spectral stabilization} \label{sec:stabilization}
The following stabilization strategy is employed throughout all subsequent simulations performed with \cref{algorithm_2}.

\begin{figure}[!ht]
	\centering
	\includegraphics[width=\linewidth]{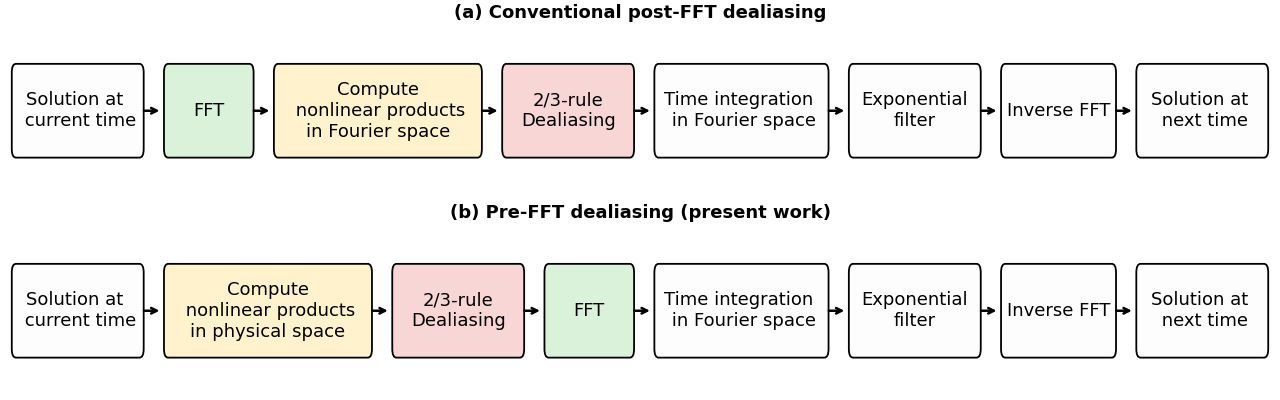}
	\caption{Comparison between the conventional post-FFT dealiasing workflow and the pre-FFT dealiasing workflow in this work.}
	\label{fig:workflow}
\end{figure}

In conventional spectral methods, the solution is first
transformed into Fourier space, where the nonlinear terms are then evaluated through operations on the Fourier coefficients. Consequently, the classical $2/3$-rule is naturally applied in Fourier space after the nonlinear terms are computed; for example, see \cite{MouhotPareschi2006}. The corresponding computational workflow is illustrated in \cref{fig:workflow}(a).

The computational workflow of \cref{algorithm_2} is fundamentally different. Instead of evaluating nonlinear terms in Fourier space, the intermediate tensors are first assembled directly in the physical domain, see \eqref{eqn:intermediate_tensors}. Only after all intermediate tensors have been constructed are they transformed into Fourier space, where their contributions are summed to evaluate the collision operator, see \eqref{Algo:15}. Because the nonlinear tensor products are already available before the multidimensional Fourier transform, the natural location of the dealiasing operation also changes. Instead of truncating Fourier modes after the transform, we apply the classical $2/3$-rule directly to the intermediate tensors prior to the FFT. The resulting workflow is shown in \cref{fig:workflow}(b). It should be emphasized that the proposed strategy does not modify the classical $2/3$-rule itself; only the stage at which it is applied is changed to match the computational structure of \cref{algorithm_2}. 

\begin{figure}[!ht]
	\centering
	\includegraphics[width=\linewidth]{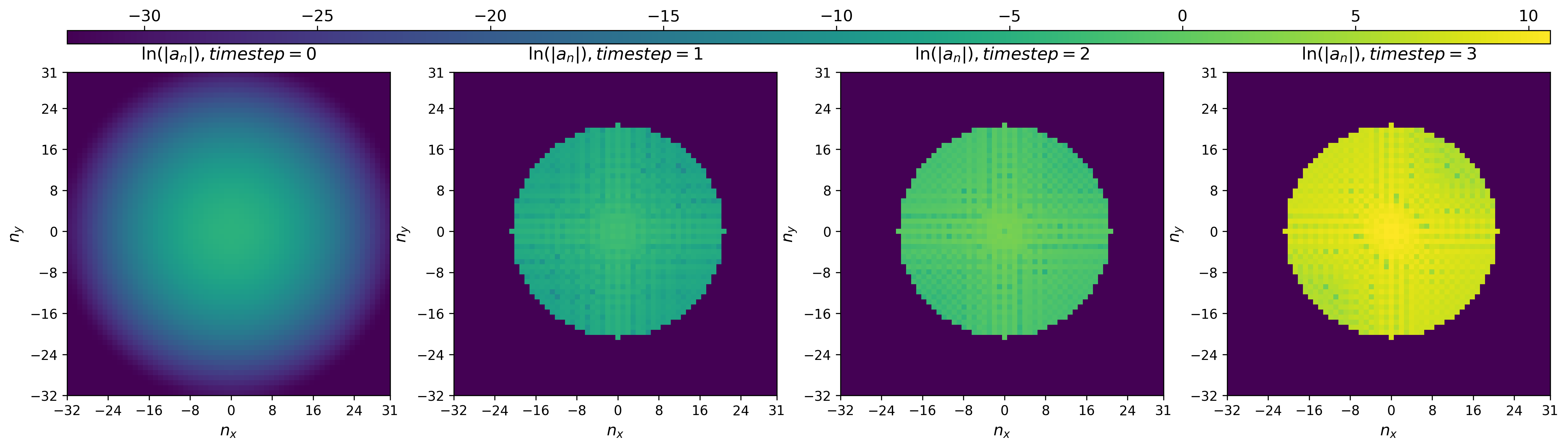}
	\caption{Evolution of the logarithmic Fourier spectrum using \cref{algorithm_2} with the $2/3$-rule dealiasing after the FFT. Spurious high-frequency modes appear after the first time step and rapidly contaminate the Fourier spectrum within only three time steps.}
	\label{fig:spectral_instability}
\end{figure}

Numerical experiments indicate that the Fourier spectra of the intermediate tensor products contain a substantial amount of high-frequency content. In the conventional workflow, these tensors are first transformed into Fourier space, and the classical $2/3$-rule is applied only after the multidimensional FFT. \cref{fig:spectral_instability} illustrates the evolution of the logarithmic Fourier spectrum obtained by using \cref{algorithm_2} with the $2/3$-rule dealiasing after the FFT. Although both the classical $2/3$-rule and an exponential spectral filter are applied after each multidimensional FFT, spurious high-frequency modes appear immediately after the first time step and rapidly propagate throughout the retained Fourier domain. Within only three time steps, the Fourier coefficients increase by several orders of magnitude over almost the entire spectrum. These results indicate that performing dealiasing only after the multidimensional FFT is insufficient to suppress the rapid growth of unresolved high-frequency modes, thereby motivating the application of the classical $2/3$-rule before the multidimensional FFT.

Note that the proposed pre-FFT dealiasing strategy is particularly suitable for problems whose solution has a compact support approximately enclosed by a ball centered at the origin of the computational domain. In this case, the intermediate tensor products inherit the same localization property, allowing the $2/3$-rule to be applied directly in the physical domain before the multidimensional FFT.

After each time integration step in Fourier space, an exponential spectral filter is applied to the retained Fourier coefficients. The filter provides an additional stabilization mechanism by smoothly damping the remaining high-frequency modes close to the cutoff while leaving the physically relevant low-frequency components nearly unchanged. Specifically,
\begin{equation}
	{a_{n}}_{\mathrm{filtered}}(k) = \sigma(k)\, a_n (k), \quad \text{where} \quad 
	\sigma(k) = \exp \left(-\eta \left( \dfrac{|k|}{k_{\max}} \right)^p \right),
\end{equation}
with $k_{\max}$ denoting the cutoff Fourier mode of the exponential filter. The parameter $\eta>0$ controls the overall damping strength, whereas the exponent $p$ determines how rapidly the filter decays near the cutoff mode. Throughout this work, we fix the filter parameters as $k_{\max} = \sqrt{2} N$, $p = 8$, and $\eta = 7.5$. These values were selected empirically based on numerical experiments to provide sufficient damping of unresolved high-frequency modes while preserving the resolved spectrum, consistent with common practice in exponential spectral filtering.

\cref{tab:parameters} summarizes the model, numerical, and stabilization parameters introduced in \cref{sec:setup} and \cref{sec:stabilization} and used in \cref{sec:qBe} and \cref{sec:3wave}.

\begin{table}[!ht]
	\centering
	\caption{Summary of the physical, numerical, and stabilization
		parameters used in the numerical simulations.}
	\label{tab:parameters}
	\begin{NiceTabular}{|l|l|c|} 
		\hline
		Parameter & Description & Fixed value \\
		\hline
		$R$ & Wave-vector domain cutoff & $R = 10$ \\
		\hline
		$M$ & Delta-function truncation parameter & $M = 4$ \\
		\hline
		$2N$ & Fourier transform sizes & $2N = 64$ in 2d case  \\  & & $2N = 32$ in 3d case\\
		\hline
		$\Delta t$ & Time step size & $\Delta t = 0.0005$ \\
		\hline
		$k_{\max}, \, p, \, \eta$ & Exponential filter parameters & $k_{\max} = \sqrt{2} N$, $p = 8$, and $\eta = 7.5$\\
		\hline
	\end{NiceTabular}
\end{table}

\subsection{Test Case 1: quantum Boltzmann equation}
\label{sec:qBe}

Unless otherwise stated, the quantum Boltzmann simulations use the
parameter values listed in \cref{tab:parameters}. We set
${\mathbf C}=1$ and $\rho=1$. The primary initial condition is
\begin{equation}
	\label{eqn:initial_power_2}
	f_0(k)=e^{-|k|^2}.
\end{equation}
To examine the influence of the decay of the initial distribution, we
also consider
\begin{equation}
	\label{eqn:initial_algebraic}
	f_0(k)=\frac{1}{1+|k|^3}.
\end{equation}

For ${\mathbf C}=1$, the collision operator contains both the quadratic
gain--loss terms
\[
f_2f_3-f_1f_2-f_1f_3
\]
and the additional linear term $-f_1$ associated with the
$1\leftrightarrow2$ quantum collision process. The numerical evolution
therefore results from a competition between the production of
excitations through the gain term and their depletion through the
quadratic and linear loss terms. Moreover, because the dispersion
relation, the collision kernel, and the initial data are all radial,
the exact equation is invariant under rotations in the wave-vector
variable. The nearly radial symmetry observed in the numerical
solutions and their Fourier spectra is therefore consistent with the
symmetry of the underlying model.

\subsubsection{Simulations in 2D}

\begin{figure}[!ht]
	\centering
	\includegraphics[width=\linewidth]{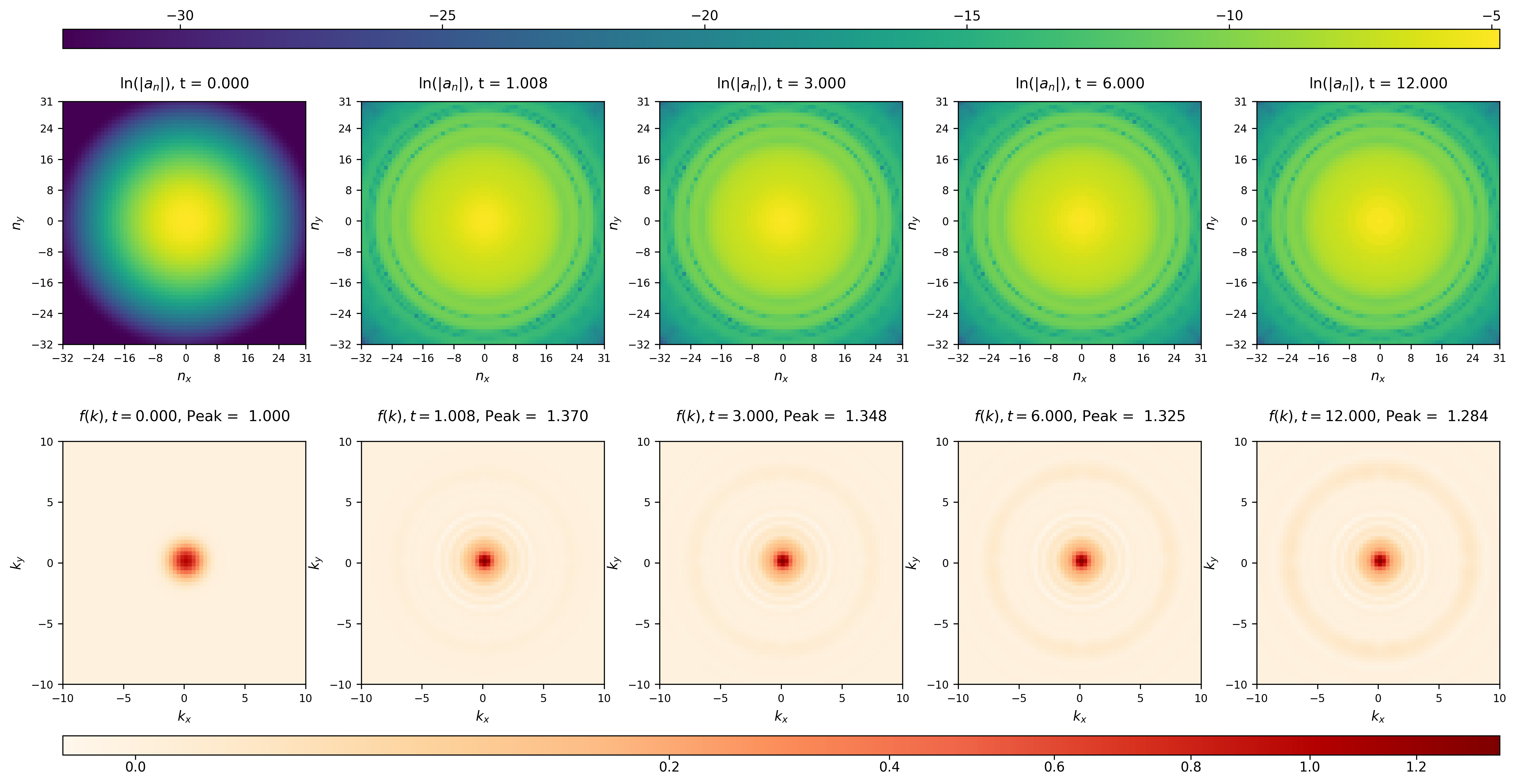}
	\caption{Physical and spectral evolution of the solution for the
		two-dimensional quantum Boltzmann equation with the polynomial
		dispersion relation \eqref{Poli1}, $\alpha=1$, and the Gaussian
		initial condition \eqref{eqn:initial_power_2}.}
	\label{fig:2D_TC1}
\end{figure}

We first consider the polynomial dispersion relation \eqref{Poli1}
with $\alpha=1$. Two initial distributions are used to illustrate how
the decay of the initial datum affects the subsequent evolution. The
results for the Gaussian initial condition
\eqref{eqn:initial_power_2} are presented in \cref{fig:2D_TC1}.

The solution remains centered at the origin and approximately radial
throughout the computation. At the beginning of the evolution, the
value at the origin increases from $1$ to approximately $1.370$ at
$t=1.008$. This transient concentration indicates that, near the
origin, the gain contribution initially dominates the corresponding
loss contributions. At later times, the maximum decreases gradually,
taking the values $1.348$, $1.325$, and $1.284$ at $t=3$, $6$, and
$12$, respectively. Thus, after the initial amplification, the
quadratic and linear loss terms progressively balance the gain term, as theoretically observed in \cite{Alonso-Gamba-Binh-2016-cauchy} (see also \cite{nguyen2017quantum} for a similar situation).

At the same time, weak concentric shells form away from the origin.
These shells indicate that the collision operator redistributes the
distribution among different wave-number scales rather than acting
only on the central mode. The logarithmic Fourier spectrum remains
radially organized but develops corresponding annular structures,
showing that the initially smooth low-frequency profile generates
additional Fourier modes during the nonlinear evolution. The absence
of substantial growth near the spectral cutoff indicates that the
stabilized algorithm remains well resolved for this simulation.

\begin{figure}[!ht]
	\centering
	\includegraphics[width=\linewidth]{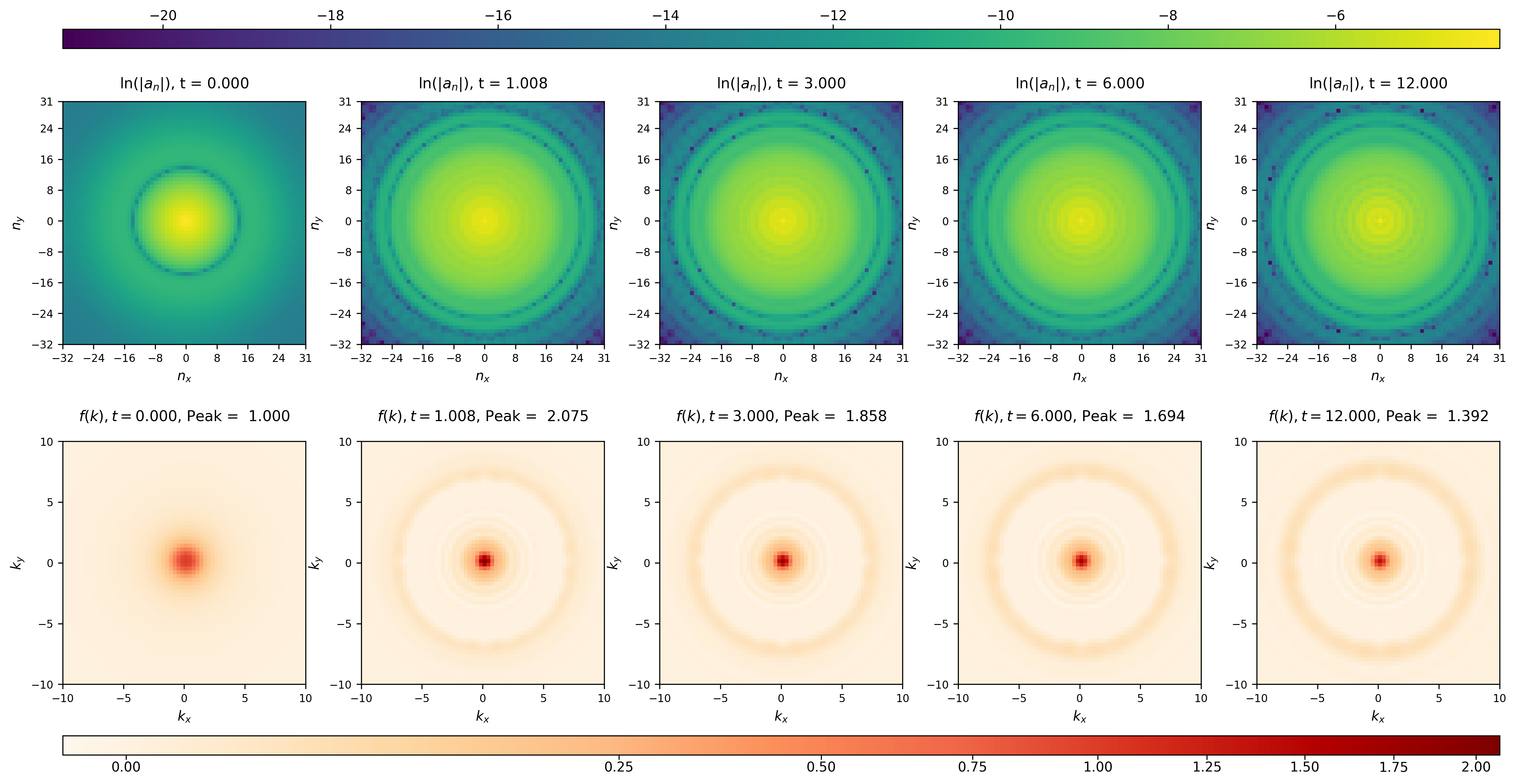}
	\caption{Physical and spectral evolution of the solution for the
		2-dimensional $C_{12}$ quantum Boltzmann equation with the polynomial
		dispersion relation \eqref{Poli1}, $\alpha=1$, and the
		algebraically decaying initial condition
		\eqref{eqn:initial_algebraic}.}
	\label{fig:2D_TC1_polynomial}
\end{figure}

The results obtained from the algebraically decaying initial condition
\eqref{eqn:initial_algebraic} are shown in
\cref{fig:2D_TC1_polynomial}.

The qualitative behavior is similar to that observed for the Gaussian
initial datum, but the transient amplification is considerably
stronger. The maximum increases from $1$ to approximately $2.075$ at
$t=1.008$ and subsequently decreases to $1.858$, $1.694$, and $1.392$
at $t=3$, $6$, and $12$, respectively. Because the algebraic initial
datum contains a larger population at moderate and large wave numbers,
more resonant configurations contribute appreciably to the collision
integral during the initial stage. This produces a stronger nonlinear
response and a more pronounced initial concentration near the origin.

The slower decay of the initial distribution also yields a broader
Fourier spectrum. Consequently, the annular structures generated
during the evolution are more visible than in
\cref{fig:2D_TC1}. Nevertheless, the solution remains approximately
radial and nonnegative, and the Fourier coefficients remain controlled
over the time interval considered. These observations indicate that
the stabilized algorithm can accommodate both rapidly and
algebraically decaying initial data in this setting.

\subsubsection{Simulations in 3D}

We next consider the three-dimensional quantum Boltzmann equation with
the polynomial dispersion relation \eqref{Poli1}, $\alpha=1$, and the
Gaussian initial condition \eqref{eqn:initial_power_2}. The displayed
panels represent 2-dimensional sections of the 3-dimensional
solution and Fourier spectrum.

\begin{figure}[!ht]
	\centering
	\includegraphics[width=\linewidth]
	{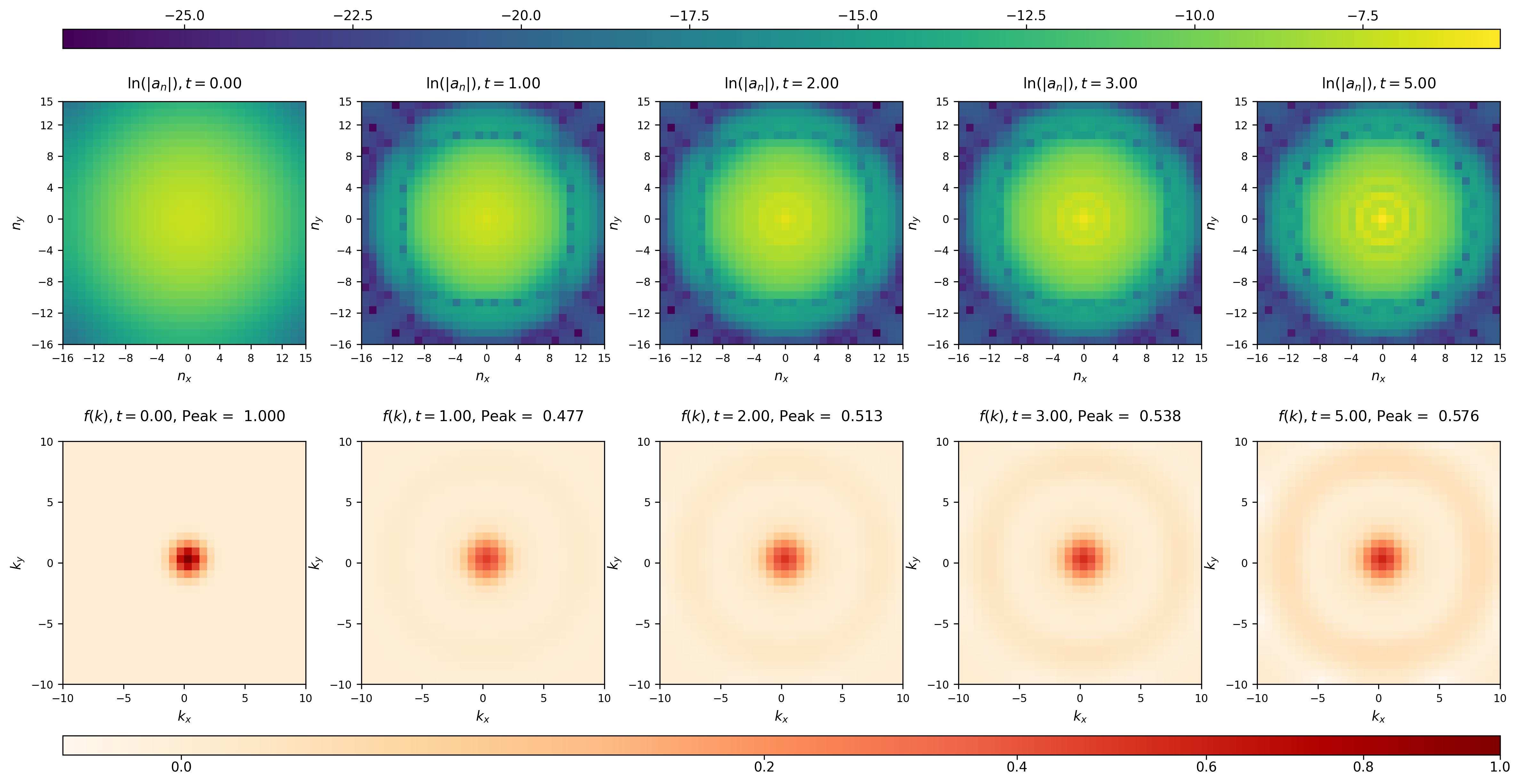}
	\caption{Physical and spectral evolution of the solution for the
		3-dimensional $C_{12}$ quantum Boltzmann equation with the polynomial
		dispersion relation \eqref{Poli1}, $\alpha=1$, and the Gaussian
		initial condition \eqref{eqn:initial_power_2}.}
	\label{fig:3D_TC1_2_80_75}
\end{figure}

As shown in \cref{fig:3D_TC1_2_80_75}, the maximum initially decreases
from $1$ to approximately $0.477$ at $t=1$. It then increases slowly
to $0.513$, $0.538$, and $0.576$ at $t=2$, $3$, and $5$,
respectively. The initial decrease is consistent with a transient
dominance of the loss contributions, while the subsequent recovery
shows that the quadratic gain term gradually repopulates the
low-wave-number region.

The solution broadens moderately and develops a weak outer shell,
while the Fourier spectrum remains concentrated around the origin and
retains an approximately radial structure. Thus, although the
dimension-dependent measure of resonant triads changes the quantitative
balance between gain and loss, the 3-dimensional computation
exhibits the same basic resonant redistribution mechanism as the
stable 2-dimensional simulations. No explosive growth of the
high-frequency Fourier modes is observed over the displayed time
interval, as theoretically observed in \cite{Alonso-Gamba-Binh-2016-cauchy} (see also \cite{nguyen2017quantum} for a similar situation). Previous numerical studies of the radial 3-wave kinetic equation provide evidence for the  energy cascade phenonmenon; see \cite{Das-Binh-2025-mumerical, Walton-Binh-2023-numerical, Walton-Binh-2024-deep, Walton-Binh-2025-numerical}. The results established here extend this picture beyond radial symmetry by numerically demonstrating energy transfer toward high frequencies in the fully non-radial setting.

\subsection{Test Case 2: 3-wave kinetic equation}
\label{sec:3wave}

The 3-wave kinetic simulations use the corresponding parameter
values reported in \cref{tab:parameters}. We set
${\mathbf C}=0$, $\rho=2$, and use the Gaussian initial condition
\eqref{eqn:initial_power_2} in all simulations.

In contrast to the quantum Boltzmann model considered in
\cref{sec:qBe}, the collision operator now contains only the quadratic
terms
\[
f_2f_3-f_1f_2-f_1f_3.
\]
Moreover, the choice $\rho=2$ causes the kernel
\[
K^{12}(k_1,k_2,k_3)
=
\left|
\omega(k_1)\omega(k_2)\omega(k_3)
\right|^2
\]
to grow rapidly with the wave numbers. The resulting semi-discrete
system is therefore strongly nonlinear and potentially very stiff.

\subsubsection{Simulations in 2D}

\begin{figure}[!ht]
	\centering
	\includegraphics[width=\linewidth]{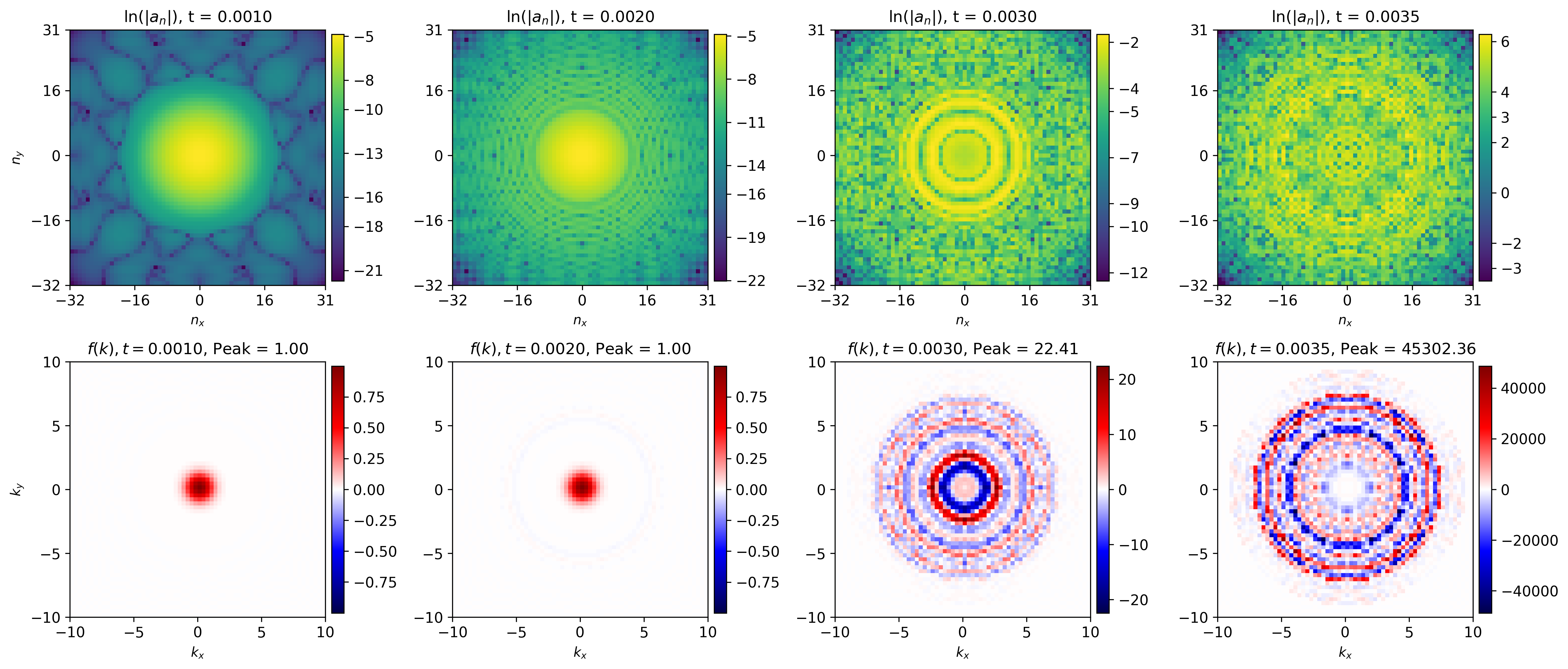}
	\caption{Physical and spectral evolution for the 2-dimensional
		3-wave kinetic equation with the polynomial dispersion relation
		\eqref{Poli1}, $\alpha=1.5$.}
	\label{fig:2D_TC2_blowup}
\end{figure}

We first consider the polynomial dispersion relation \eqref{Poli1}
with $\alpha=1.5$. The corresponding results are presented in
\cref{fig:2D_TC2_blowup}.

During the first two displayed time levels, the solution remains close
to the initial Gaussian profile and its maximum remains approximately
equal to $1$. By $t=0.003$, a sequence of alternating positive and
negative concentric shells has formed, and the maximum has increased
to approximately $22.41$. At $t=0.0035$, the maximum reaches
approximately $4.53\times10^4$. Simultaneously, the initially localized
Fourier spectrum spreads across almost the entire retained frequency
domain and increases by several orders of magnitude.

The rapid transfer of amplitude toward higher Fourier modes is
consistent with the strong wave-number dependence of the kernel when
$\rho=2$ and $\omega(k)=|k|^{1.5}$.  The computation therefore exhibits an apparent blow-up, which may be interpreted as a transfer of energy (energy cascade) toward increasingly high frequencies in the solution; see \cite{Soffer-Binh-2019-energy,staffilani2025finite} for related theoretical results.

\begin{figure}[!ht]
	\centering
	\includegraphics[width=\linewidth]{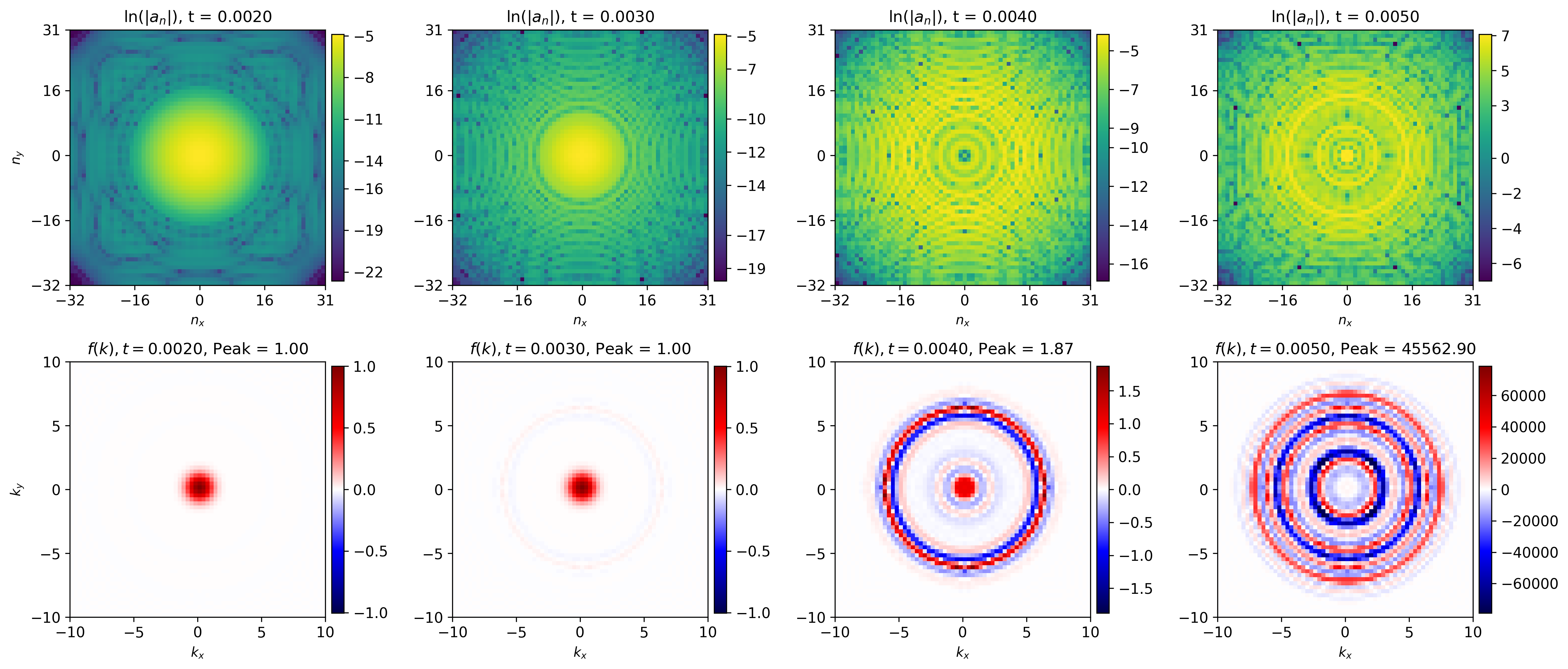}
	\caption{Physical and spectral evolution for the 2-dimensional
		3-wave kinetic equation with the mixed polynomial dispersion
		relation \eqref{Poli2}, where $\alpha=1$, $\beta=1.5$,
		$c_1=0.6$, and $c_2=0.4$.}
	\label{fig:2D_TC2_poli2_form}
\end{figure}

We next consider the mixed polynomial dispersion relation
\eqref{Poli2} with
\[
\alpha=1,\qquad
\beta=1.5,\qquad
c_1=0.6,\qquad
c_2=0.4.
\]
The corresponding results are shown in
\cref{fig:2D_TC2_poli2_form}.

The mixed dispersion relation changes the geometry of the resonance
manifold and delays the onset of the rapid growth, but it does not
eliminate it. The solution remains close to its initial profile at
$t=0.002$ and $t=0.003$. At $t=0.004$, a pronounced oscillatory shell
appears at an intermediate wave number and the maximum increases to
approximately $1.87$. By $t=0.005$, several alternating shells are
present and the maximum reaches approximately
$4.56\times10^4$. The Fourier spectrum likewise develops concentric
structures and becomes broadband.

The linear component $c_1|k|$ and the superlinear component
$c_2|k|^{1.5}$ generate a resonance geometry different from that of
the homogeneous dispersion relation used in
\cref{fig:2D_TC2_blowup}.   Thus, the comparison between
\cref{fig:2D_TC2_blowup,fig:2D_TC2_poli2_form} demonstrates that the
dispersion relation strongly affects the transient resonant pattern.

\subsubsection{Simulations in 3D}

Finally, we consider the 3-dimensional 3-wave kinetic equation
with the polynomial dispersion relation \eqref{Poli1},
$\alpha=1.5$. The displayed panels are 2-dimensional sections of the
3-dimensional numerical solution and its Fourier spectrum.

\begin{figure}[!ht]
	\centering
	\includegraphics[width=\linewidth]{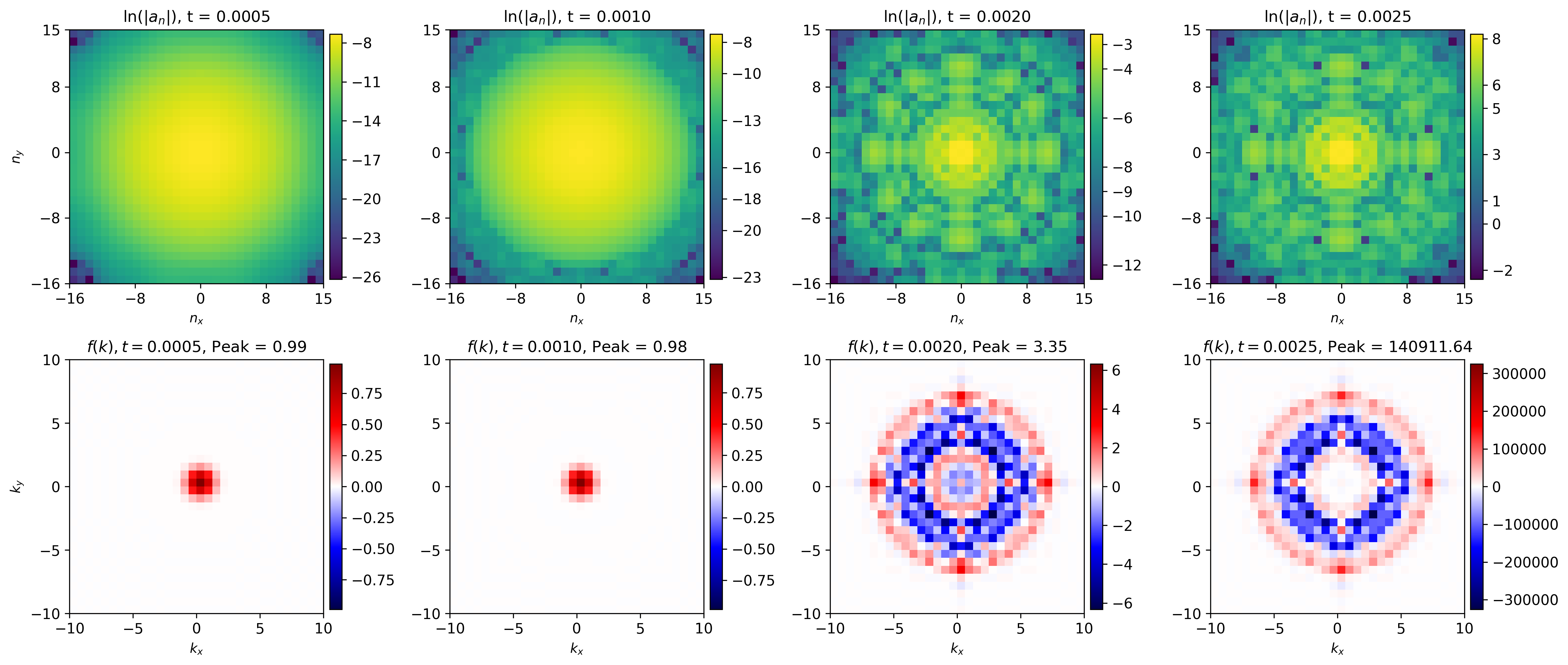}
	\caption{Physical and spectral evolution for the 3-dimensional
		3-wave kinetic equation with the polynomial dispersion relation
		\eqref{Poli1}, $\alpha=1.5$.}
	\label{fig:3D_TC2_blowup}
\end{figure}

The 3-dimensional computation exhibits the same qualitative
instability as the 2-dimensional case, but on an even shorter time
scale. The maximum remains close to $1$ at $t=0.0005$ and
$t=0.001$. At $t=0.002$,   the
maximum has increased to approximately $3.35$. By $t=0.0025$, the
maximum reaches approximately $1.41\times10^5$, and the solution is
dominated by large oscillations distributed over an annular region.
The Fourier spectrum simultaneously loses its initial localization
and develops substantial amplitude throughout the retained high modes.

The stronger and earlier growth is consistent with the increased
number of interacting configurations available in three dimensions,
together with the rapidly growing kernel corresponding to
$\rho=2$.   Consequently, \cref{fig:3D_TC2_blowup} may be interpreted as indicating a transfer of energy toward high frequencies (an energy cascade) in the solution; see \cite{Soffer-Binh-2019-energy,staffilani2025finite} for related theoretical results. Previous numerical studies of the radial 3-wave kinetic equation have provided evidence of an energy cascade; see \cite{Das-Binh-2025-mumerical, Walton-Binh-2023-numerical, Walton-Binh-2024-deep, Walton-Binh-2025-numerical}. The numerics developed here shows that this cascade mechanism remains present in the fully non-radial numerical schemes.

\section{Conclusion}
\label{sec:conclusion}

In this work, we have developed spectral algorithms for multidimensional 3-wave kinetic equations and \(C_{12}\) quantum Boltzmann equations with general polynomial dispersion relations. The resonance constraint was treated through a truncated Fourier representation of the Dirac distribution supported on the resonance manifold. This formulation incorporates the collision operator into a Fourier spectral framework without requiring an explicit parametrization or a piecewise-polynomial approximation of the resonance manifold.

We derived 2  numerical methods. The first is a direct spectral algorithm in which the distribution function is approximated by a truncated Fourier series and the evolution equations for its Fourier coefficients are evaluated through precomputed kernel coefficients. Its overall computational complexity is
\[
\mathcal{O}\bigl(L(2N)^{3d}\bigr).
\]
The second algorithm follows a different computational principle: multidimensional Fourier transforms are applied directly to the kernel-weighted nonlinear terms assembled in the physical variables. This formulation reduces the computational complexity to
\[
\mathcal{O}\bigl(L(2N)^{2d}\log(2N)\bigr).
\]
The numerical comparison in 2 dimensions shows that the 2 methods produce solutions agreeing to nearly machine precision. At the same time, the fast spectral algorithm yields substantial computational savings, with measured speedups increasing from approximately \(93\) for \(2N=16\) to more than \(2200\) for \(2N=64\). The fast method also remains computationally feasible for the three-dimensional resolutions considered in this work.

The direct construction of the kernel-weighted nonlinear tensors generates significant high-frequency content and therefore requires an appropriate stabilization procedure. We introduced a pre-FFT implementation of the classical \(2/3\)-rule, adapted to the computational structure of the fast algorithm, together with an exponential filter applied to the retained Fourier coefficients after each time step. The numerical experiments demonstrate that applying dealiasing only after the multidimensional Fourier transform is insufficient to prevent the rapid contamination of the retained spectrum. By contrast, the proposed combination of pre-FFT dealiasing and exponential spectral filtering substantially improves the stability of the computations for distributions localized near the center of the computational domain.

For the \(C_{12}\) quantum Boltzmann equation, the stabilized algorithm captures the competition between the quadratic gain term and the quadratic and linear loss terms. The numerical solutions retain the approximate radial symmetry dictated by the radial dispersion relations, collision kernels, and initial data. In 2 dimensions, both Gaussian and algebraically decaying initial distributions exhibit an initial amplification near the origin followed by a gradual decrease of the maximum. The algebraically decaying initial datum produces a stronger transient response and a broader Fourier spectrum because it contains a larger population at moderate and high wave numbers. In 3 dimensions, the solution initially decreases near the origin and subsequently recovers as the gain mechanism repopulates the low-wave-number region. In all these simulations, the Fourier spectrum remains controlled over the time intervals considered. 

The 3-wave kinetic equation displays a markedly different behavior. For the rapidly growing collision kernel corresponding to \(\rho=2\), the computations develop alternating oscillatory shells, rapid growth of the solution amplitude, and pronounced broadening of the Fourier spectrum. These features appear in both 2 and 3 dimensions and may be interpreted as numerical evidence of an apparent transfer of energy toward high frequencies, consistent with an energy-cascade mechanism. The 3-dimensional computation exhibits this growth on a shorter time scale than the corresponding 2-dimensional computation. Moreover, the comparison between homogeneous and mixed polynomial dispersion relations shows that the geometry of the resonance manifold strongly affects the transient resonant patterns and may delay, although not necessarily suppress, the onset of rapid high-frequency growth. Our results confirm that the energy-cascade mechanism previously observed numerically under radial symmetry also occurs in the non-radial 3-wave kinetic equation.

Several questions remain for future investigation. A rigorous consistency and convergence analysis of the approximation of the resonance constraint would provide a mathematical foundation for the proposed algorithms. It would also be useful to develop positivity-preserving and conservative discretizations, adaptive choices of the truncation parameter \(M\), and higher-order or implicit time-integration methods for the stiff 3-wave regime. Moreover, further reductions in memory consumption will be important for higher-resolution 3-dimensional simulations. 

\section*{Acknowledgements}
This research is partly funded by University of Economics Ho Chi Minh City (UEH), Vietnam. T. T. Le acknowledges this financial support.  M.-B. T. is  funded in part by      NSF CAREER  DMS-2303146, and NSF Grants DMS-2204795, DMS-2305523,  DMS-2306379.

\bibliographystyle{plain}
\bibliography{References_new}

\end{document}